\def\reals{{{\rm l}\kern - .15em {\rm R} }} 

\documentclass[12pt]{article}
\usepackage{amssymb,amsmath}

\newtheorem{theorem}{Theorem}[section]
\newtheorem{lemma}[theorem]{Lemma}
\newtheorem{proposition}[theorem]{Proposition}

\newenvironment{proof}[1][Proof]{\begin{trivlist}
\item[\hskip \labelsep {\bfseries
#1}]}{\end{trivlist}}
\newenvironment{definition}[1][Definition]{\begin{trivlist}
\item[\hskip \labelsep {\bfseries
#1}]}{\end{trivlist}}

 \newcommand{\qed}{\nobreak
\ifvmode \relax \else
\ifdim\lastskip<1.5em \hskip-\lastskip
\hskip1.5em plus0em minus0.5em \fi \nobreak \vrule height0.75em width0.5em
depth0.25em\fi}

\def\reals{{{\rm l}\kern - .15em {\rm R} }} 

\begin{document}

\begin{centering}

{\huge BRST EXTENSION OF }

{\huge \vspace{.5cm} 
GEOMETRIC QUANTIZATION}
\\

\renewcommand{\thefootnote}{\fnsymbol{footnote}}
\vspace{1.4cm}
{\large Ronald Fulp}\\

\vspace{.5cm}

 Department of Mathematics, North Carolina State University, Raleigh,
NC 27695-8205.\\
\it 
E-mail: fulp@math.ncsu.edu\\
\rm
\vspace{.5cm}

\rm
\vspace{.5cm}
\vspace{.5cm}

\begin{abstract}

 Consider a physical system for which a mathematically rigorous geometric quantization procedure exists. Now subject the system to a finite set of irreducible first class (bosonic) constraints. It is shown that there is a mathematically rigorous BRST quantization of the constrained system whose cohomology at ghost number zero recovers the constrained quantum states. Moreover this space of constrained states has a well-defined Hilbert space structure inherited from that of the original system. Treatments of these ideas in the Physics literature are more general but suffer from having states with infinite
or zero  "norms" and thus are not admissible as states. Also the BRST operator for many systems require regularization to be well-defined. In our more restricted context we show that our treatment does not suffer from any of these difficulties.

\end{abstract}
\end{centering}
\newpage
\section{Introduction}
Let $M$ be a symplectic manifold with symplectic form $\omega$ and for $f,g\in C^{\infty}M$ let $[f,g]$ denote the induced Poisson brackets of $f$ and $g.$ Assume that $\sigma$ is a nonnegative polarization on $M$ with real directions (see Woodhouse\cite{W},  pages 184-186). Let $Q=M/\sigma$ and let ${\cal H_{\sigma}}$ denote the Hilbert space of $\sigma$-wave functions $\tilde s =s\nu$ where $s$ is a section of the prequantum bundle $B$ over $M$ and $\nu$ is a section of an appropriate square-root line bundle $\delta_{\sigma}$ over $M.$ Sections of these bundles must respect the polarization described in detail in Woodhouse in order to obtain a bundle over $Q=M/{\sigma}.$ The space  ${\cal H}_{\sigma}$ is a space of sections of the bundle $B_{\sigma}$ over $Q$ defined by $B_{\sigma}=B\otimes \delta_{\sigma}.$ Woodhouse denotes such sections by $\tilde s$ but we generically denote them  by $\psi$ and refer to them as {\bf  {${\sigma}$-wave functions}} or simply as {\bf wave functions}. 

Let $C_{\sigma}^{\infty}M$ denote the linear space of smooth real-valued functions $f$ such that the flow of the Hamiltonian vector field  of $f$ preserves the polarization $\sigma.$ It is not difficult to see that this space is closed under Poisson brackets. Moreover, a Dirac-like correspondence exists from the space $C_{\sigma}^{\infty}M$ to self-adjoint operators on ${\cal H}_{\sigma}$  relative to a natural inner-product on 
${\cal H_{\sigma}}.$  If $f\in C_{\sigma}^{\infty}M$  we denote the corresponding operator on  $C_{\sigma}^{\infty}M$ by $\hat f$ and observe that $\widehat{[f,g]}=i[\hat f,\hat g]$ for such $f$ and $g.$

In the present work we address the problem of quantizing this system when the symplectic manifold $M$  is subjected to additional constraints. Assume that one has a set of constraint functions $\{G_a\}  \quad (a=1,2,\cdots, m)$ defined on $M$ subject to the following conditions: \newline

\noindent (1) the constraints are first-class, i.e., there exist smooth functions $\{C_{ab}^c\}$ on $M$ having the property that $\{G_a,G_b\}=C_{ab}^cG_c,$
\newline
(2) $\Sigma=\{p\in M | \quad G_a(p)= 0, a=1,2,\cdots,m\}$ is a submanifold of $M$ such that at each point of $M$ there is an open subset $U$ of $M$ on which the $\{G_a \}$ define a chart on $U\cap \Sigma$ and are also the first $m$ coordinates of a chart of $M$ defined on $U,$
\newline
(3) the constraints are irreducible, and
\newline
(4) the flow of the Hamiltonian vector field $X_a$ of each constraint $G_a$ preserves the polarization $\sigma.$
\newline 
Thus the constraints are irreducible, regular,  and define a self-adjoint operator $\hat G_a$ on the Hilbert space of wave-functions $\psi.$ 
\bigskip

\bigskip

To quantize this constrained system we employ the methods of BRST cohomology.
Thus we must develop the BRST machinery required to quantize the constrained system which we do by following Henneaux and Teitelboim \cite{HT}.
In particular, we initially follow the development on page 319 which we now outline briefly for the convenience of the reader. They first extend the  phase space $M$  to include new variables the so-called ghosts $\{\eta^a \}$ which are new anti-commuting variables assumed to be in  one-one correspondence with the constraints. Additionally they require corresponding ghost momenta $\{ {\cal P}_a \}.$ Roughly, their new phase space is $M\oplus <\eta^a,{\cal P}_b>$ where $<\eta^a,{\cal P}_b>$ is the complex linear space spanned by $\eta^a,{\cal P}_b,$ but we will see that this description is not fully adequate to describe the situation.  The observables on this extended phase space are complex polynomials in the variables 
$\{\eta^a\}$ and $\{ {\cal P}_a\}$ with coefficients in $C^{\infty}M.$ Henneaux and Teitelboim denote this space of observables by $C(\eta^a)\otimes C^{\infty}M\otimes C ({\cal P}_a).$ In fact this space is a graded algebra. The ghosts and their momenta are assigned an odd parity (see page 190 \cite{HT}) while the elements of $C^{\infty}M$ are even. Essentially, then, they have a polynomial algebra in which elements of $C^{\infty}M$ commute with every polynomial while the $\eta^a$ and ${\cal P}_b,$ at the classical level,  satisfy the relation
$$ \eta^a{\cal P}_b+{\cal P}_b\eta^a=-\delta ^a_b.$$
The ghost momenta are regarded as pure imaginary supernumbers in a superalgebra whereas the ghosts are real supernumbers.
Under quantization $\hat \eta^a$ and $\hat {\cal P}_a$ may be identified as  operators defined on an enlarged space of wave functions. In fact this enlarged space is not a Hilbert space as it possesses a  scalar product which  is degenerate. It is our intent to rigorously  describe these extended wave functions, the relevant scalar product, and the operators corresponding to the classical observables. 
\bigskip

\section{ Extended Phase Space.} Let $\Lambda$ denote the superalgebra of supernumbers modeled on a Grassmann algebra generated  by either a finite or countably infinite number of generators. Our conventions regarding supernumbers and their properties subscribe to those of Rogers \cite{R}. This algebra has a $Z_2$ grading $\Lambda=\Lambda^0 \oplus \Lambda^1$ where, as usual,  $x\in \Lambda^0$ is assigned parity $\varepsilon(x)=0$ while $\eta\in \Lambda^1$ is assigned parity $\varepsilon(\eta)=1.$ Generally a Lagrangian in this context is a mapping from the tangent bundle of some configuration supermanifold ${\cal M}$ into $\Lambda.$ Thus in case ${\cal M}=\Lambda^{r|s}=(\Lambda^0)^r \times (\Lambda^1)^s,$ $L$ is  locally a function of the coordinates $(q^i,\theta^{\alpha},\dot {q^i},  \dot{\theta^{\alpha}})$ and generically is required to be  even and real (maps into $\Lambda^0_{Re}$, see \cite{HT}) so that the ``momenta"  $p_i=\frac{\partial L}{\partial q^i}$ are even and real while $\pi_{\alpha}=\frac{\partial L}{\partial \theta^{\alpha}}$ are imaginary and odd. A  corresponding phase space is then constructed along with observables  which are functions of $(q^i,\theta^{\alpha},p_i,\pi_{\alpha})$ as is usual in the Hamiltonian formalism.

For our purposes the  {\bf extended configuration space} will be $(M/{\sigma})\times (\Lambda_{Re}^1)^m$ where $m$ is the number of constraints  $\{G_a\}.$ We identify our so-called {\bf extended phase space} with the space $M\times (\Lambda_{Re}^1)^m\times ( \hat \Lambda_{Im}^1)^m$ where $\hat \Lambda^1$ denotes a formal dual of $\Lambda^1,$ i.e., $\hat \Lambda^1$ is a copy of $\Lambda^1$ where elements of $\Lambda^1$ are denoted $\eta^a$ while those of $\hat \Lambda^1$
are denoted ${\cal P}_b$ and are regarded as the momenta of $\eta^a.$ A contraction $\eta^a{\cal P}_a$ is simply an element of $\Lambda.$ {\bf Classical observables} for us are mappings $F$ from extended  phase space $M\times (\Lambda_{Re}^1)^m\times ( \hat \Lambda_{Im}^1)^m$  into a tensor product $\Lambda \hat \otimes B_{\sigma} \hat \otimes \hat \Lambda$ such that each mapping $F$ is a finite sum of functions $F_{rs}$ homogenous in the $\eta^a$'s and ${\cal P}_b$'s:
$$F_{rs}(x,\eta^a,{\cal P}_b)=\eta^{a_1}\eta^{a_2}\cdots\eta^{a_r} 
f^{a_1a_2\cdots a_r}_{b_1b_2\cdots b_s}(x){\cal P}_{b_1}{\cal P}_{b_2}\cdots{\cal P}_{b_s}$$
where $f^{a_1a_2\cdots a_r}_{b_1b_2\cdots b_s}$ is in $C_{\sigma}^{\infty}M.$  Recall that $C^{\infty}_{\sigma}M$ denotes the linear space of all smooth real-valued functions $f$ on $M$ whose Hamiltonian vector fields preserve the polarization ${\sigma}.$
The range of these classical observables is the tensor product indicated above but the tensor product  is subject to symmetries defined by requiring that the following commutation relations hold: 
$$[\eta^a,f]=\eta^af-f\eta_a=0 , [{\cal P}_a,f]={\cal P}_a f-f{\cal P}_a=0$$
and 
$$[{\cal P}_a,\eta^b]={\cal P}_a \eta^b+ \eta^b {\cal P}_a= [\eta^b,{\cal P}_a]=-\delta^b_a$$
where $f \in C_{\sigma}^{\infty}M.$

Here and hereafter we adhere to the summation convention. Additionally, we often use a multi-index notation so that $F_{rs}=\eta^Af^B_A{\cal P}_B$ where $A=(a_1,a_2,\cdots, a_r), B=(b_1,b_2,\cdots, b_s),$ and 
$$\eta^A=\eta^{a_1}\eta^{a_2}\cdots\eta^{a_r} \quad \quad 
{\cal P}_B={\cal P}_{b_1}{\cal P}_{b_2}\cdots {\cal P}_{b_s}.$$
In this notation, the multi-indices $A$ and $B$ are increasing multi-indices of positive integers and there is an implied sum over the individual indices of each multi-index. When we wish to consider observables of the form $F=F_{rs}=\eta^Af^B_A{\cal P}_B$ where $A=(a_1,a_2,\cdots, a_r), B=(b_1,b_2,\cdots, b_s),$ 
we will say that  $F$ is {\bf homogeneous}  and  in such a case we say that {\bf $F$ has degree
$(r,s)$ or  $r+s$} depending on the context. Often we denote the observables defined by 
$$(x,\eta^a,{\cal P}_b)\rightarrow \eta^a, \quad  \quad(x,\eta^a,{\cal P}_b)\rightarrow {\cal P}_b, \quad \quad (x,\eta^a,{\cal P}_b)\rightarrow f^B_A(x)$$
simply by $\eta^a,{\cal P}_b,$ and $f^B_A,$ respectively.

The set of all classical observables  ${\cal O}$ is a linear space over the real numbers, moreover there is a well-known Poisson bracket defined on this space (see \cite{HT} page 146). The bracket satisfies the graded Jacobi identity and additionally satisfies the conditions:
\begin{eqnarray}
&[F,G]=-(-1)^{\varepsilon_F \varepsilon_G}[G,F] \\ & 
[FG,H]=F[G,H]+(-1)^{\varepsilon _G\varepsilon_H}[F,H]G \\ & 
\varepsilon[F,G]=\varepsilon F+\varepsilon G
 \\
& 
[F,G]^*=-[G^*,F^*]
 \end{eqnarray}
for all $F,G,H\in {\cal O}.$ Here conjugation is required to be linear and to satisfy the conditions $(\eta^a)^*=\eta^a, ({\cal P}_b)^*=-{\cal P}_b, $ and $(zw)^*=w^*z^*.$ In particular, 
$$(\eta^Af^B_A{\cal P}_B)^*=   ({\cal P}_B)^* (f^B_A)^* (\eta^A)^*$$
where, for example, $(\eta^A)^*=\eta^{a_r}\eta^{a_{r-1}}\cdots \eta^{a_1}.$

\bigskip

\section {Extension of  Geometric Quantization.} The space of {\bf extended geometric quantum states} will be denoted throughout the paper by ${\cal S}.$ A function $\psi$ denotes such a state iff it is a mapping from our configuration space ${\cal M}=(M/{\sigma})\times (\Lambda^1_{Re})^m$ into $B_{\sigma}\hat \otimes \Lambda$ such that for $(x,\eta^a)\in {\cal M}$
$$\psi(x,\eta^a)=\psi_0(x)+\psi_a(x)\eta^a+\psi_{ab}(x)\eta^a\eta^b+\cdots + \psi_{12\cdots m}(x)(\eta^1\eta^2 \cdots \eta^m).$$
In multi-index notation we write $\psi=\psi_I (x)\eta^I$
where the sum extends only over increasing multi-indices. Here for each increasing multi-index I,
$\psi_I$ denotes a section of the bundle $B_{\sigma}$ defined by  Woodhouse\cite{W} (pages 185,186). Recall from the introduction that, following Woodhouse, we call sections of the bundle $B_{\sigma}\rightarrow M/{\sigma},$ {\bf  ${\sigma}$-wave functions}. 
Woodhouse shows that these ${\sigma}$-wave functions form a pre-Hilbert space via an inner product defined by identifying the inner product of two such ${\sigma}$-wave functions $s\nu$ and $s'\nu'$ as an integral 
$$\int_{M/{\sigma}}(s,s')\nu\nu'$$ where $(s,s')$ is the inner product of sections of the Hermitean prequantization bundle $B$ over $M.$  He does this by showing that $(s,s')\nu\nu'$ can be identified as a top form on $M.$ 

If $\psi, \phi \in {\cal S},$ then for multi-indices $I,J$ the components $\phi_I,\psi_J$ are ${\sigma}$-wave functions and we denote by  $(\psi_I | \phi_J)$ their inner product as defined above by Woodhouse. We will abuse Woodhouse's notation by also writing $(\psi|\phi)=(\psi_I\eta^I|\phi_J\eta^J)=(\psi_I|\phi_J)\eta^I\eta^J.$ We reserve the notation $(\psi,\phi)$ for a {\it scalar product} which we define below on our space ${\cal S}$ of states. This scalar product is degenerate and will eventually be utilized to obtain an inner product on the appropriate space of BRST cohomology classes. 

For $\psi,\phi\in {\cal S},$ {\it define} the {\bf scalar product} of $\psi$ and $ \phi$ 
by $$(\psi,\phi)=top\{(\psi_I |\phi_J)(\eta^I)^*\eta^J\}.$$ Here {\it top} denotes the coefficient of $\eta^1\eta^2 \cdots \eta^m$ the term of highest degree in $(\psi| \phi).$ In a less condensed notation $(\psi,\phi)$ is the sum of all terms of the form 
$$top\{ (\psi_{a_1a_2\cdots a_r} |\phi_{b_1b_2\cdots b_s})
(\eta^{a_1}\eta^{a_2}\cdots\eta^{a_r})^*\eta^{b_1}\eta^{b_2}\cdots\eta^{ b_s}\}$$
where $r+s=m.$

Now that we have a space of states our next objective is to create quantum observables. Thus we must assign to each $F\in{\cal O}$ a linear operator on the space ${\cal S.}$

Woodhouse has shown that if $f\in C_{\sigma}^{\infty}M$ is real-valued, then there is a well-defined operator $\hat f$
(he uses the notation $\tilde f$) on the completion of the space of the {\it ${\sigma}$-wave functions} which is self-adjoint relative to the inner product $(\cdot | \cdot )$ defined on his Hilbert space 
${\cal H}_{\sigma}.$ We extend this 
operator to ${\cal S}$  in the obvious way 
$$\hat f(\psi_J\eta^J)=\hat f(\psi_J)\eta^J$$
and note that it is self-adjoint (see the formal definition below) relative to the scalar product: 
$$(\hat f(\psi),\phi)=top\{ (\hat f(\psi_J) | \phi_I )(\eta^J)^*\eta^I\}= top\{ (\psi_J |\hat f(\phi_I ))(\eta^J)^*\eta^I\}
=(\psi,\hat f \phi).$$
 
The operators corresponding to the special classical observables $\eta^a$ and ${\cal P}_b$ are defined by  
$$\hat \eta^a\psi=\eta^a\psi \quad \quad \hat {\cal P}_a\psi=
\frac {1}{i}\frac {\partial \psi}{\partial \eta^a}.$$ 

One can show by a direct calculation that $\hat \eta^a \hat \eta^b=-\hat \eta^b\hat\eta^a$ and that
$\hat {\cal P}_a\hat {\cal P}_b=-\hat {\cal P}_b\hat {\cal P}_a;$ consequently,
$$[\hat \eta^a, \hat \eta^b]=0  \quad \quad                    [ \hat {\cal P}_a,\hat {\cal P}_b]=0.$$

Also note that for $\psi\in {\cal S},$
$$(\hat {\cal P}_a\hat \eta^b)(\psi)=\frac{1}{i}\frac {\partial}{\partial \eta^a}(\eta^b\psi)=
\frac {1}{i}\delta^b_a\psi-\eta^b\hat {\cal P}_a(\psi)$$ 
and 
$$[\hat {\cal P}_a,\hat \eta^b]=-i\delta^b_a.$$

The general prescription for quantizing classical observables $F,G\in {\cal O}$ requires that
$$\widehat {[F,G]}=i[\hat F,\hat G].$$
 We see that the commutators we have derived above satisfy this condition. Moreover if we require that for $f\in C_{\sigma}^{\infty}M$ that 
$$[\hat f,\hat \eta^a]=0=[\hat f, \hat {\cal P}_b],$$ 
then the required condition holds for these classical observables.

To establish the general case we require that if $F=\eta^A f_A^B{\cal P}_B,$ then 
$\hat F=\hat {\eta^A} \hat {f_A^B}\hat {{\cal P}_B},$ i.e., we require that the $\hat {\cal P}_b$ operators
act first then  operators of the form $\hat f$ for $f\in C_{\sigma}^{\infty}M$ followed by the action of the  operators $\hat \eta^a.$ Here, for each multi-index $A=(a_1,a_2,\cdots ,a_r),$
$$\hat \eta^A=\hat\eta^{a_1}\hat\eta^{a_2}\cdots\hat \eta^{a_r} 
 \quad \quad\hat {\cal  P}_A=\hat{\cal  P}_{a_1}\hat{\cal  P}_{a_2}\cdots \hat{\cal  P}_{a_r} .$$

With these conventions we show that the required quantum relation holds for all classical observables $F,G.$  To accomplish this and for other purposes as well we first show that $\widehat {FG}=\hat F\hat G.$ Assume first that $F=\eta^A f_A^B{\cal P}_B, G=\eta^I g_I^J{\cal P}_J$
are homogeneous, i.e., that the implied sum over the  multi-indices A,B,I,J extend only over multi-indices of a fixed degrees, thus the sums extend over
$$A=(a_1,a_2,\cdots ,a_r),B=(b_1,b_2,\cdots,b_s), I=(i_1,i_2,\cdots ,i_u),J=(j_1,j_2,\cdots ,j_v)$$
for fixed values of $r,s,u,v.$ 
Notice that to get $FG=\eta^A f_A^B{\cal P}_B\eta^I g_I^J{\cal P}_J$ in the correct order to compute $\widehat {FG}$ one must rewrite the terms ${\cal P}_B\eta^I$ in reversed order using the commutator relation $[{\cal P}_b,\eta^i]=-\delta^i_b$ repeatedly. If we define $\tilde \eta^i=\hat \eta^i$  and $\tilde {\cal P}_b=(-i)\hat {\cal P}_b,$ then these new variables satisfy the same commutator relations as the original classical variables $\eta^i, {\cal P}_b.$ Consequently, we may write $\widetilde{{\cal P}_B\eta^I}=\tilde{\cal P}_B\tilde \eta^I$ where a tilde over
a quantity with a multi-index means that each factor is replaced with the tilde of that factor without changing the order of the factors, for example,  $\tilde{\cal P}_B=\tilde{\cal  P}_{b_1}\tilde{\cal  P}_{b_2}\cdots \tilde{\cal  P}_{b_s} .$ It easily follows that $\widetilde{FG}=\tilde F \tilde G.$ On the other hand 
$$(-i)^{s+v}\widehat {FG}=\widetilde{FG}=\tilde F \tilde G=(-i)^{s+v}\hat F \hat G$$
and so $\widehat{FG}=\hat F \hat G $  for homogeneous $F,G.$
Now in general $FG$ is a sum of homogeneous terms and since $H\rightarrow \hat H$ is a linear mapping it follows that  $\widehat{FG}=\hat F \hat G $ holds for all $f,G\in {\cal O}.$
This shows that the mapping from ${\cal O}$ into linear operators on ${\cal S}$ is a homomorphism of associative algebras.

It is well-known (see \cite{HT} pages 146,235) that the Poisson bracket on ${\cal O}$ satisfies the identity
$$[FG,H]=F[G,H]+(-1)^{\varepsilon_G \varepsilon_H}[F,H]G,$$
Moreover, this same identity holds for ``graded commutators" of operators; consequently,
 $$[\hat F \hat G,\hat H]=\hat F[\hat G,\hat H]+(-1)^{\varepsilon_G \varepsilon_H}[\hat F,\hat H]\hat G$$
 (also recall that, by definition, $\varepsilon(\hat K)=\varepsilon (K)$ for every $K\in {\cal O}$).

Now observe that the identity $\widehat {[F,G]}=i[\hat F,\hat G]$ holds for the ``generators"
${\cal P}_b,\eta^a, f\in C_{\sigma}^{\infty}M$ of each classical observable $H\in {\cal O}.$
To see that the quantization identity $\widehat {[F,G]}=i[\hat F,\hat G]$ holds for arbitrary classical observables observe
first that it suffices to show this for homogenous $F,G \in {\cal O}.$ We indicate why the identity holds for the homogeneous case  via 
an inductive argument. If $H=\eta^{a_1}\eta^{a_2}\cdots\eta^{a_r}H_{a_1a_2\cdots a_r}^{b_1b_2\cdots b_s}{\cal P}_{b_1}{\cal P}_{b_2}\cdots {\cal P}_{b_s}\in {\cal O}$ is homogeneous recall that the {\bf degree of }$H$ is   $deg(H)=r+s.$  Note that the quantization identity holds for all $F,G$ such that $deg(F)+deg(G)=1$ since in that case either $F$  or $G$ is in
$C_{\sigma}^{\infty}M.$ Now assume that the required identity holds for all $F,G$ such that 
$deg(F)+deg(G)<k+1$ and let $F,G$ be observables such that $deg(F)+deg(G)=k+1.$ Since $[F,G]=\pm[G,F]$ it is no loss of generality to assume that $deg(F)>1.$ Write $F=F_1F_2$ with 
$deg(F_i)+deg)(G)<k+1, i=1,2.$ Then,
$$\widehat {[F,G]}=\widehat {[F_1F_2,G]}=\widehat {F_1[F_2,G]}(-1)^{\varepsilon_{F_2}\varepsilon_G}\widehat{[F_1,G]F_2}$$$$=\hat F_1 i[\hat F_2,\hat G]+(-1)^{\varepsilon_{F_2}\varepsilon_G} i[\hat F_1,\hat G]\hat F_2=i[\hat F_1\hat F_2,\hat G]=i[\hat F,\hat G]$$

It follows that the quantization condition holds for every classical observable.

\section{Adjoint Operators}

In this section we consider how adjoints of operators on quantum state space ${\cal S}$ are defined. Recall that our scalar product on ${\cal S}$ is highly degenerate so care must be taken to assure that the notion is well-defined. For this purpose we need the following lemma.

\begin{lemma}
Assume that $\psi$ is a state in ${\cal S}$ such that $(\psi,\phi)=0$ for every state $\phi.$
Then $\psi=0.$
\end{lemma}

\begin{proof} Let $\psi=\psi_I \eta^I$ and consider one specific multi-index $I_0$ of the sum. Let $J_0$
denote the complement of $I_0$ in the sequence $(1,2,\cdots, m)$ (recall that $m$ is the number of constraints). Consider an arbitrary homogeneous element $\phi=\phi_{J_0}\eta^{J_0}$ of  ${\cal S}$ and notice that $0=(\psi,\phi)=(\psi_{I_0}|\phi_{J_0})\eta^{I_0}\eta^{J_0}.$ Consequently, 
$(\psi_{I_0}|\phi_{J_0})=0$ for every ${\sigma}$-wave function $\phi_{J_0}$ and since $(\quad|\quad)$ is positive definite it follows that $\psi_{I_0}=0.$ Since the multi-index $I_0$ was arbitrary, it follows that $\psi=0.$ 
\end{proof}

\begin{definition} If $A:{\cal S}\rightarrow {\cal S}$ is a linear mapping on state space and if there exists a linear operator $B$ from ${\cal S}$ to ${\cal S}$ such that $(A(\psi),\phi)=(\psi,B(\phi))$ for all states
$\psi,\phi \in {\cal S},$ then we say 
$B$ is a {\bf right-adjoint} of $A.$ Similarly we say $B$ is a {\bf left-adjoint}  of $A$ if 
$(B(\psi),\phi)=(\psi,A(\phi))$ for all states
$\psi,\phi \in {\cal S}.$ We say that $B$ is an {\bf adjoint } of $A$ if it is both a left and right adjoint. It follows from the lemma that if left or right adjoints exist, then they are unique.  If $A$ has an adjoint it is denoted by  $A^{\dagger}.$  An operator $A$ is {\bf self-adjoint} iff $A^{\dagger}=A$
and {\bf skew-adjoint} iff $A^{\dagger}=-A$
\end{definition}

\noindent{\bf  Remark.} Notice that if $A,B$ are linear operators on ${\cal S}$ and if they have  adjoints, then so does the composite operator $AB$ and $(AB)^{\dagger}=B^{\dagger}A^{\dagger}.$ The analogous formulae hold for left and right adjoints.
 
We now show that $\hat \eta^a$ is self-adjoint and that $\hat {\cal P}$ is skew-adjoint with respect to the scalar product.

\begin{proposition} The operator $\hat \eta^a$ is self-adjoint and $\hat {\cal P}_r$ is skew-adjoint relative to the scalar
product $(\quad,\quad)$ on ${\cal S}$ for each $1\leq a,r \leq m.$
\end{proposition}

\begin{proof} Let $\psi,\phi\in {\cal S}$ and note that 
$$ (\hat \eta^a(\psi),\phi)=(\psi_I\eta^a\eta^I,\phi_J\eta^J)=top\{(\psi_I|\phi_J)
(\eta^a\eta^I)^*\eta^J\}$$$$=top\{(\psi_I|\phi_J)
(\eta^I)^*\eta^a\eta^J\}
=(\psi,\hat \eta^a(\phi)). $$
Consequently, $\hat \eta^a$ is self-adjoint. Next we show $\hat {\cal P}_r$ is skew-adjoint.
We find it useful in this particular proof to modify the notation we have been using for states $\psi$
by writing $\psi=\psi_I\eta^I$ where the multi-indices $I$ are permitted to vary over {\it all skew-symmetric} multi-indices. Thus the coefficients of $\psi$ and $\phi$ below are not the same as in the 
rest of the paper since they include a factor $\frac{1}{r!}$ where $r$ is the number of components of the relevant multi-index. With this in mind, observe that 
$$top(\psi|\hat {\cal P}_r(\phi))=\frac{1}{i} top\{(\psi |(\phi_r+2\phi_{ra}\eta^a+3\phi_{rbc}\eta^b\eta^c+\cdots))\}$$$$=
\frac{1}{i} top(\sum_{s+t=m+1}\{t(\psi_{a_1a_2\cdots a_s}|\phi_{rb_1\cdots b_{t-1}})
(\eta^{a_s}\eta^{a_{s-1}}\cdots\eta^{a_1}\eta^{b_1}\eta^{b_2}\cdots \eta^{b_{t-1}}))$$
and that  \newline $top(\hat {\cal P}_r(\psi)|\phi)=$
$$
(-\frac{1}{i})top(\{\sum_{s+t=m+1}s(\psi_{ra_1a_2\cdots a_{s-1}} |\phi_{b_1b_2\cdots b_t})
(\eta^{a_{s-1}}\cdots\eta^{a_2}\eta^{a_1}\eta^{b_1}\eta^{b_2}\cdots \eta^{b_t}))\}.
$$
On the other hand since $s+t>m$
$$0=\hat {\cal P}_r\{\sum_{s+t=m+1}(-1)^{s-1}(\psi_{a_1a_2\cdots a_s}|\phi_{b_1\cdots b_t})
(\eta^{a_s}\cdots\eta^{a_2}\eta^{a_1}\eta^{b_1}\eta^{b_2}\cdots \eta^{b_t})\}$$
$$=\frac {1}{i}\{\sum_{s+t=m+1}(-1)^{s-1}s(\psi_{a_1a_2\cdots a_{s-1}r}|\phi_{b_1\cdots b_t})
(\eta^{a_{s-1}}\cdots\eta^{a_2}\eta^{a_1}\eta^{b_1}\eta^{b_2}\cdots \eta^{b_t})$$
$$+ \sum_{s+t=m+1}(-1)^{s-1}(-1)^s(\psi_{a_1a_2\cdots a_s}|t\phi_{rb_1\cdots b_{t-1}})
(\eta^{a_s}\cdots\eta^{a_2}\eta^{a_1}\eta^{b_1}\eta^{b_2}\cdots \eta^{b_{t-1}})\}$$
$$=\frac {1}{i}\{\sum_{s+t=m+1}(-1)^{s-1}(-1)^{s-1}s(\psi_{ra_1a_2\cdots a_{s-1}}|\phi_{b_1\cdots b_t})
(\eta^{a_{s-1}}\cdots\eta^{a_2}\eta^{a_1}\eta^{b_1}\eta^{b_2}\cdots \eta^{b_t})$$
$$- \sum_{s+t=m+1}(-1)^{s-1}(-1)^{s-1}t(\psi_{a_1a_2\cdots a_s}|\phi_{rb_1\cdots b_{t-1}})
(\eta^{a_s}\cdots\eta^{a_2}\eta^{a_1}\eta^{b_1}\eta^{b_2}\cdots \eta^{b_{t-1}})\}.$$
This implies that
$0=-top(\hat {\cal P}_r(\psi)|\phi)-top(\psi|\hat {\cal P}_r(\phi))$ and   that $0=top(\hat {\cal P}_r(\psi)|\phi)=-top(\psi|\hat {\cal P}_r(\phi).$ Consequently,  
$(\hat {\cal P}_r(\psi),\phi)=(\psi,-\hat {\cal P}_r(\phi))$ and  $\hat {\cal P}_r$ is skew-adjoint.

\end{proof}

We now return to our usual convention that multi-indices are always increasing finite sequences of
positive integers.

\begin{proposition} If  $F=\eta^A f_A^B{\cal P}_B \in {\cal O}$ is a classical observable, then
$\hat F$ has an adjoint and 
$$(\hat F)^{\dagger}=(-1)^{|B|}(\hat{\cal P}_B)^{\dagger}\hat f_A^B(\hat \eta^A)^{\dagger}=(-1)^{|B|}(-1)^{{[\frac {|A|}{2}}][\frac{|B|}{2}]} \hat {\cal P}_B \hat f_A^B\hat \eta^A$$ 
where $|A|,|B|$ are the number of terms in the respective multi-indices and $[x]$  denotes the greatest integer in the number $x.$
\end{proposition}

\begin{proof} The proposition is an immediate consequence of the last proposition and the remark preceding it. 
\end{proof}

\section{Ghost number and BRST Operators}

The goal of the BRST program is to obtain the constrained quantum states  as BRST cohomology classes. One must first construct two operators on the space of states subject to certain conditions. First one needs a mapping $\Omega$ from the space of states to itself such that $\Omega^2=0.$ This provides a differential for the theory so that the required cohomology is simply $$H(\Omega)=kernel ( \Omega)/  boundary (\Omega).$$ Additionally one needs a ghost number operator ${\cal G}.$ This operator provides a grading on the cohomology complex $H(\Omega).$ Elements of  $H(\Omega)$ of ghost number $g$ are denoted $H^g(\Omega).$  For a successful encoding of the physics it  turns out that the quantum states {\it must} reside precisely at ghost number zero cohomology, $H^0(\Omega).$ In addition to these requirements $\Omega$ must be self-adjoint and ${\cal G}$ must be skew-adjoint (at least this is one formula for success). These conditions restrict not only the operators $\Omega$ and ${\cal G}$ but they also restrict the kind of scalar product available. It can be shown (see \cite{HT} page 299) that the scalar product must be degenerate if all of these conditions are to be true. In fact if $\psi$ has ghost number $p$ and $\psi'$ has ghost number $q,$ then the scalar product of $\psi$ and $\psi'$ must be zero unless $p+q=0$ and it is this fact which requires consideration of a scalar product of the type we have constructed (following \cite{HT}) above.

Following Henneaux and Teitelboim once again (page 298, bosonic irreducible case) we define the {\bf ghost number operator}  ${\cal G}$ by:
$${\cal G}=\frac {i}{2}[\hat \eta^a\hat {\cal P}_a-\hat {\cal P}_a\hat \eta^a].$$
In physics terminology the {\bf  $\eta^a$'s are called ghosts} and the ghost number operator essentially keeps track of the number of ghosts present in a given quantum state. 

\begin{proposition} The operator ${\cal G}$ is skew-adjoint. Moreover if \newline $\psi=\psi_{a_1a_2\cdots a_s}\eta^{a_1}\eta^{a_2} \cdots \eta^{a_s}\in {\cal S}$
is homogeneous, then 
$${\cal G}\psi=-\frac{m}{2}+s.$$
\end{proposition} 
\begin{proof} First notice that ${\cal G}$ is skew-adjoint since 
$${\cal G}^{\dagger}=(-\frac {i}{2})[\hat{\cal P}_a^{\dagger}(\hat \eta^a)^{\dagger}-(\hat \eta^a)^{\dagger}\hat {\cal P}_a^{\dagger}]=(-\frac{i}{2}[\hat \eta^a\hat {\cal P}_a-\hat {\cal P}_a\hat \eta^a]=-{\cal G}.$$
Next observe that \newline ${\cal G}\psi_0=\frac{i}{2}[-\hat {\cal P}_a\hat \eta^a](\psi_0)=$
$$(-\frac{1}{2})[\frac {\partial}{\partial \eta^1}(\eta^1\psi_0)+\frac {\partial}{\partial \eta^2}(\eta^2\psi_0)+\cdots+\frac {\partial}{\partial \eta^m}(\eta^m \psi_0)]=-\frac{m}{2}\psi_0.$$
Finally, note that since 
$${\cal G}=\frac {i}{2}[\hat \eta^a\hat {\cal P}_a-\hat {\cal P}_a\hat \eta^a]=\frac {i}{2}[2\hat \eta^a\hat {\cal P}_a+i\delta^a_a]=-\frac{m}{2}+i\hat \eta^a\hat {\cal P}_a$$
and $ \eta^a\frac{\partial}{\partial \eta^a}(\eta^{a_1}\eta^{a_2}\cdots\eta^{a_r})=
r (\eta^{a_1}\eta^{a_2}\cdots\eta^{a_r}),$ it follows that 
$${\cal G}(\psi_{I_r}\eta^{I_r})=\eta^a\frac{\partial}{\partial \eta^a}(\psi_{I_r}\eta^{I_r})-\frac{m}{2}(\psi_{I_r}\eta^{I_r})$$
and ${\cal G}(\psi)=(r-\frac {m}{2})\psi.$ The proposition follows.

\end{proof}

Since our goal is to realize the constrained quantum states as BRST cohomology classes we must understand how the classical BRST operator relates to the classical BRST observable $\Omega.$
We do not repeat the construction of the classical BRST operator in detail as this is done rigorously on pages 194-195 and 224-226 of the book by Henneaux and Teitelboim \cite{HT}. We do need to know that the classical observable $\Omega$ is constructed inductively in such a manner that 
$\Omega=\Omega^{(0)}+\Omega^{(1)}+\cdots$ where $\Omega^{(0)}=\eta^aG_a,$  
$$\Omega^{(p)}=\eta^{b_1}\eta^{b_2}\cdots\eta^{b_{p+1}}f^{a_1a_2\cdots a_p}_{b_1b_2\cdots b_{p+1}}
{\cal P}_{a_p}{\cal P}_{a_{p-1}}\cdots {\cal P}_{a_1},$$ and where in the inductive process, we include
the requirement that the $f^{a_1a_2\cdots a_p}_{b_1b_2\cdots b_{p+1}}$ be real.
In physpeak the ${\cal P}_b$'s are called {\bf anti-ghosts so one says that $\Omega^{(p)}$ has
anti-ghost number} p.

In the inductive construction of the $\Omega^{(p)}$ one assumes that 
$\Omega^{(1)},\Omega^{(2)},\cdots, \Omega^{(p)}$ have been constructed and then one constructs $\Omega^{(p+1)}.$ This involves the Koszul-Tate operator $\delta.$ Proper consideration of this operator would lead us too far afield so we refer to Henneaux and Teitelboim \cite{HT} for details regarding $\delta.$
To find $\Omega^{(p+1)}$ once $\Omega^{(1)},\Omega^{(2)},\cdots, \Omega^{(p)}$  are known one defines an observable $D^{(p)}$ by
$$D^{(p)}=\frac{1}{2}\{\sum_{k=0}^p[\Omega^{(p-k)},\Omega^{(k)}]_{C_{\sigma}^{\infty}M}
+\sum_{k=0}^{p-1}[\Omega^{(p-k+1)},\Omega^{(k)}]_{{\cal P},\eta}\}$$
where $[\quad,\quad]_{C_{\sigma}^{\infty}M}$ denotes Poisson brackets induced from $M$ ignoring the ghosts and anti-ghosts and $[\quad, \quad]_{{\cal P},\eta}$ denotes Poisson brackets involving the ghosts and anti-ghosts and ignoring the brackets of elements of $C_{\sigma}^{\infty}M.$ It is then shown that there exists an observable $\Omega^{(p+1)}$ such that $\delta(\Omega^{(p+1)})=-D^{(p)}$ and that {\it any such solution will suffice} to construct the full BRST operator  having the required property that 
$\Omega^2=0.$

Since it is essential to us that the quantum operator $\hat \Omega$ not only have square zero but that it be self-adjoint, we show this latter property in detail.

\begin{lemma} If $F\in{\cal O}, F=\eta^Af^B_A{\cal P}_B,$  $(f^B_A)^*=f^B_A,$ and $F^*=F,$ then $\hat F^{\dagger}=\hat F.$
\end{lemma}
\begin{proof}
It suffices to consider homogeneous $F\in{\cal O},$
$$F=\eta^{a_1}\eta^{a_2}\cdots\eta^{a_r}f_A^B{\cal P}_{b_1}{\cal P}_{b_2}\cdots {\cal P}_{b_s}.$$
Then 
$$F=F^*=(-1)^s{\cal P}_{b_s}{\cal P}_{b_{s-1}}\cdots {\cal P}_{b_1}f_A^B\eta^{a_r}\eta^{a_{r-1}}\cdots\eta^{a_1},$$
since $(\eta^a)^*=\eta^a$ and $({\cal P}_b)^*=-{\cal P}_b.$ But  $(\hat\eta^a)^*=\hat\eta^a$ and $(\hat{\cal P}_b)^*=-\hat{\cal P}_b$ so
$$\hat F ^{\dagger}=(\hat\eta^{a_1}\hat\eta^{a_2}\cdots \hat\eta^{a_r} \hat f_A^B \hat{\cal P}_{b_1}\hat{\cal P}_{b_2}\cdots \hat{\cal P}_{b_s})^{\dagger}= \quad \quad\quad\quad\quad\quad\quad\quad\quad $$$$\quad\quad\quad\quad\quad(-1)^s\hat{\cal P}_{b_s}\hat{\cal P}_{b_{s-1}}\cdots \hat{\cal P}_{b_1} \hat f_A^B\hat\eta^{a_r}\hat\eta^{a_{r-1}}\cdots\hat\eta^{a_1}=\widehat{F^*}=\hat F.$$
\end{proof}

\begin{proposition} If ${\Omega}$ is the classical BRST operator and if the homological perturbation of
$\Omega$ is finite, $\Omega=\Omega^{(0)}+\Omega^{(1)}+\cdots \Omega^{(n)},$ then $\Omega$ is a self-adjoint quantum observable.
\end{proposition}
\begin{proof} We first show that the classical observable $\Omega$ is real by an inductive argument.
Note first that $(\Omega^{(0)})^*=(\eta^a G_a)^*=G_a^*(\eta^a)^*=\eta^aG_a=\Omega^{(0)}$ and $\Omega^{(0)}$ is real. Assume that $\Omega^{(1)},\Omega^{(2)},\cdots, \Omega^{(p)}$ are real. We show that $\Omega^{p+1}$ is real. Observe that 
$$(D^{(p)})^*=\frac{1}{2}\{\sum_{k=0}^p[\Omega^{(p-k)},\Omega^{(k)}]^*_{C_{\sigma}^{\infty}M}
+\sum_{k=0}^{p-1}[\Omega^{(p-k+1)},\Omega^{(k)}]^*_{{\cal P},\eta}\}.$$
Since $$[F^*,G^*]=-[G,F]^*=(-1)(-1)(-1)^{\varepsilon_F\varepsilon_G}[F,G]^*$$ and $\varepsilon(\Omega^{(k)})=1,$ we see that ${D^{(p)}}^*=-D^{(p)}.$
Now Henneaux and Teitelboim show that there exists $\tilde\Omega^{(p+1)}$ such that 
$\delta(\tilde\Omega^{(p+1)})=-D^{(p)}$ and that any other solution of this latter equation suffices. If 
$\tilde\Omega^{(p+1)}$ is not real, define $\Omega^{(p+1)}=\frac{1}{2}[\tilde\Omega^{(p+1)}+(\tilde\Omega^{(p+1)})^*],$ then $\Omega^{(p+1)}$ is real and $\delta(\Omega^{(p+1)})=\frac{1}{2}[\delta(\tilde\Omega^{(p+1)})+\delta(\tilde\Omega^{(p+1)})]=
\frac{1}{2}[-D^{(p)}-D^{(p)}]=-D^{(p)}$ as required. By induction we see that each summand  $\Omega^{(p)}$ of $\Omega$ is real. It follows from the lemma that 
$\hat \Omega=\hat \Omega^{(0)}+\hat \Omega^{(1)}+\cdots \hat \Omega^{(n)}$ is self-adjoint and
$[\hat\Omega,\hat\Omega]=\widehat{[\Omega,\Omega]}=0.$  Since $\hat\Omega$ is odd, $\hat\Omega^2=0.$

\end{proof}

\bigskip

It is shown in Henneaux and Teitelboim under general conditions that if one has a scalar product on the space ${\cal S}$ of states and if $\psi,\phi$ are states having ghost numbers $p$ and $q,$ respectively, then $(\psi,\phi)$ can be nonzero only when $p+q=0.$ The aim of BRST cohomology in this instance is to produce physical states at ghost number zero. If 
$$\psi=\psi_0+\psi_a\eta^a+\psi_{ab}\eta^a\eta^b+\cdots$$
is a state, then $\psi_0$ plays the role of a traditional Dirac state at least when $M=T^*Q$ for some configuration manifold $Q.$ Now ${\cal G}(\psi_0)=-\frac {m}{2},$ so to obtain states with ghost number zero one first shows that  the states at ghost number $-\frac{m}{2}$ are dual, with respect to the scalar product, to those at ghost number $\frac{m}{2}.$ If, for each nonnegative integer $k,$ we let ${\cal S}^k$ denote the linear space of all states $\psi\in {\cal S}$ such that $\psi=\psi_{a_1a_1\cdots a_p}\eta^{a_1}\eta^{a_2}\cdots \eta^{a_p}, $ then the relevant states reside in ${\cal S}^0\oplus {\cal S}^m$ (notice the two different gradings). To obtain states at ghost number zero one extends the space of classical observables and its corresponding space of quantum observables along with the requisite BRST machinery to obtain a new BRST operator $\Omega+\tilde\Omega$ such that 

$$H^0(\Omega+\tilde\Omega)=[H^{-\frac {m}{2}}(\Omega) \otimes H^{\frac {m}{2}}(\tilde\Omega) ]                                                                                       
      \oplus [H^{\frac {m}{2}}(\Omega) \otimes H^{-\frac {m}{2}}(\tilde\Omega)].$$
      
 One then shows that  $H^0(\Omega+\tilde\Omega)$ is indeed the space of constrained quantum states.  This pathway works quite generally, but many obstacles can occur. The operator $\Omega$ itself may require regularization. The scalar product $(\psi,\psi)$ may be infinite and require regularization. 
 
 It is our contention that these problems are avoided if one applies the BRST construction to a system which already admits a rigorous geometric quantization but is then subjected to further first class constraints. 
 
 To finish our program, we must first show that ${\cal S}^0$ and ${\cal S}^m$ are dual and to prove this we first need to set the stage with a few observations.    
      
For each $p$ define a pairing of  ${\cal S}^p$ and ${\cal S}^{m-p}$ by
mapping an ordered pair $\{\psi,\phi\}\in  {\cal S}^p \times {\cal S}^{m-p}$ to the scalar product 
$(\psi,\phi)$ of the two states. Notice that the {\it proof} of Lemma 4.1 shows that the mapping from ${\cal S}^p$ to the dual of ${\cal S}^{m-p}$ defined by $\psi \rightarrow \alpha_{\psi},$ where $\alpha_{\psi}(\phi)=(\psi,\phi),$ is an injective mapping. In case $p=0,$ clearly ${\cal S}^p$ can be identified with the Hilbert space ${\cal P}{\cal W}$ of ${\sigma}$-wave functions. Similarly for  $p=m$ we have that 
$${\cal S}^m=\{\psi\in{\cal S} | \psi=
\psi_{12\cdots m}(\eta^1\eta^2\cdots \eta^m), \psi_{12\cdots m}\in {\cal P}{\cal W} \}$$
can be identified with the space  ${\cal P}{\cal W}$ of ${\sigma}$-wave functions by simply dropping the term $\eta^1\eta^2\cdots \eta^m.$ To avoid confusion let $\zeta$ be the inverse of this identification, specifically $\zeta$ is the mapping from ${\cal S}^m$ to ${\cal S}^0={\cal P}{\cal W}$ defined by $\zeta(\psi)=\psi_{12\cdots m},$ for $\psi\in {\cal S}^m.$ Notice that for $\psi_0\in {\cal S}^0$ and $\phi \in {\cal S}^m,$ we have $(\phi,\psi_0)=(\zeta(\phi) | \psi_0).$ Consequently, the scalar product of an element of ${\cal S}^0$ and an element of ${\cal S}^m$ is essentially the inner product in the Hilbert space ${\cal P}{\cal W}$ of ${\sigma}$-wave functions. Using these facts we can prove the following theorem.

\begin{theorem} The BRST cohomology at ghost number $-\frac{m}{2}$ is isomorphic to the BRST cohomology at ghost number $\frac {m}{2}.$ Both cohomology spaces are isomorphic to the Hilbert space ${\cal P}{\cal W}$ of ${\sigma}$-wave functions.
\end{theorem}
\begin{proof}
First observe that since ${\cal S}^0$ is identified with ${\cal P}{\cal W},$ and since 
$\zeta(\Omega({\cal S}^{m-1}))$ is also a subspace of the same space we have 
$${\cal S}^0=\zeta(\Omega({\cal S}^{m-1}))\oplus \zeta(\Omega({\cal S}^{m-1}))^{\bot}$$ where the orthogonal decomposition is taken relative to the inner product on  ${\cal P}{\cal W}.$ Moreover, notice that the BRST cohomology at ghost number $-\frac{m}{2}$ is simply the kernel of the restriction of $\Omega$ to ${\cal S}^0$ while the BRST cohomology at ghost number $\frac {m}{2}$ is ${\cal S}^m/\Omega({\cal S}^{m-1})$ which may be identified with the dual space $ (\zeta(\Omega({\cal S}^{m-1}))^{\bot})^*.$ Define a mapping 
$\Lambda: ker(\Omega|_{{\cal S}^0)})\longrightarrow  (\zeta(\Omega({\cal S}^{m-1}))^{\bot})^*$ by $\Lambda_{\psi_0}(\phi)=(\psi_0|\phi).$ Notice first that if $\Lambda_{\psi_0}=0,$ then $(\psi_0|\phi)=0$ for all $\phi \in \zeta(\Omega({\cal S}^{m-1}))^{\bot}.$ But it is also true that $(\psi_0|\rho)=(\psi_0 | \zeta(\Omega(\phi)))=( \psi_0,\Omega(\phi))=(\Omega(\psi_0),\phi)=0$ for all $\rho=\zeta(\Omega(\phi))\in \zeta(\Omega({\cal S}^{m-1}),$ since $\Omega$ is self-adjoint. Consequently, $\psi_0$ is an element of the Hilbert space ${\cal P}{\cal W}$ which is orthogonal to every other element and  so is zero. It follows that $\Lambda$ is injective. We show that it is also surjective. Let
$h\in  (\zeta(\Omega({\cal S}^{m-1}))^{\bot})^*$ so that $h$ is a continuous linear mapping from $ \zeta(\Omega({\cal S}^{m-1}))^{\bot}$ into the complex numbers ${\bf C}.$ Define $\tilde h: {\cal P}{\cal W}\rightarrow {\bf C}$ by $\tilde h(x)=h(x)$ for $x\in (\zeta(\Omega({\cal S}^{m-1}))^{\bot})$ and define $\tilde h(x)=0$ for $x\in \zeta(\Omega({\cal S}^{m-1})).$ Clearly this uniquely defines an element of the dual of  ${\cal P}{\cal W}.$ Thus there exists $\psi_o\in {\cal P}{\cal W}={\cal S}^0$ such that $\tilde h(\phi)=(\psi_0|\phi)$ for all $\phi \in {\cal S}^0.$ Thus $(\psi_0, \Omega({\cal S}^{m-1}))=(\psi_0 |\zeta(\Omega({\cal S}^{m-1})))=\tilde h(\zeta(\Omega({\cal S}^{m-1})))=0.$ But since $\Omega$ is self-adjoint we have that $(\Omega(\psi_0),{\cal S}^{m-1})=(\psi_0,\Omega({\cal S}^{m-1})=0.$ By the proof of Lemma 4.1, it follows that $\Omega(\psi_0)=0$ and that $\psi_0\in ker (\Omega |_{{\cal S}^0}).$ It now follows that $\Lambda_{\psi_0}(\phi)=(\psi_0|\phi)=\tilde h(\phi)=h(\phi),$ for all $\phi\in \zeta(\Omega({\cal S}^{m-1}))^{\bot}.$ It follows that $\Lambda(\psi_0)=h$ and $\Lambda$ is surjective.
The theorem follows.

\end{proof}

\section{The Constrained Quantum States} 
At this point we enlarge the space of classical observables and the BRST machinery to obtain the goal outlined just a couple of paragraphs prior to the statement of the last theorem.

Our extended manifold is $M\times T^*S^m$ where $S^m$ denotes the $m$-sphere (recall that $m$ is the number of constraints). The manifold $M\times T^*S^m$ is given its natural product structure obtained from that of $M$ and $T^*S^m.$ This manifold clearly admits a polarization whose leaves are $L\times T_q^*S^m$ where $L$ is a leaf of the polarization of $M$ and $T_q^*S^m$ is the fiber of $T^*S^m$ over $q\in S^m.$ One has a Poisson bracket defined on $C^{\infty}(M\times T^*S^m)$ with respect to the extended structure. Functions $f\in C^{\infty}(M\times T^*S^m)$ which are constant on leaves of $M\times T^*S^m$ will be denoted by $C^{\infty}((M/{\sigma})\times S^m)$ where $(M/{\sigma})\times S^m$ plays the role of our new ``configuration" space. Geometric quantization now applies to these functions to produce operators $\hat f$ on the enlarged space of ${\sigma}$-wave functions which we denote by $\widetilde {\cal PW}.$

If $(q^i)$ denotes local coordinates on $M/{\sigma}$ and $Q^a$ denotes local coordinates on $S^m,$ then we denote their conjugate momenta coordinates by $(p_j)$ and $(P_b)$ respectively. The new constraints for the combined system are:
$$G_a=0 \quad\quad\quad P_b=0.$$
Thus the $P_b$'s play a double role being both constraints and momenta. We introduce new ghosts and ghost momenta $\tilde \eta^a, \tilde {\cal P}_b$ corresponding to the new constraints $P_b=0.$ In the physics literature the ghost momenta have negative ghost number, are usually denoted by $\overline C_b$ and are called anti-ghosts but we will not need this language for our purposes (although it is clearly useful when applying the BRST formalism). In a local notation our extended phase space has ``coordinates"
$$(q^i, Q^a, p_j, P_b,\eta^c, \tilde \eta^d, {\cal P}_e, \tilde {\cal P}_f).$$
Classical extended observables are functions from 
$$M\times T^*S^m\times (\Lambda_{Re}^1)^m\times ({\hat\Lambda_{Im}^1})^m$$ to $B_{\sigma} \hat\otimes \Lambda$ where $B_{\sigma}$ is the new line bundle obtained from geometric quantization of the extended system. We do {\it not} require all such 
functions, however, just certain ones  of the form:
$$F(x,y,\eta^c, \tilde \eta^d, {\cal P}_e, \tilde {\cal P}_f)=\eta^Af_A^B(x){\cal P}_B\otimes \tilde\eta^Ig_I^J(y)\tilde{\cal P}_J.$$
where $f^B_A\in C_{\sigma}^{\infty}M$ and $g_I^J$ is in the corresponding subspace of $C^{\infty}T^*S^m.$

The space of states of $T^*S^m$ will be denoted by $\tilde{\cal S}$ and consists of all maps from $S^m\times (\Lambda_{Re}^1)^m$ into $\Lambda$ defined by 
$$\tilde \psi =\tilde \psi_0+\tilde\psi_a\tilde\eta^a +\cdots + 
\tilde \psi_{123\cdots m}(\tilde\eta^1\tilde\eta^2\cdots \tilde\eta^m)$$
where $\tilde\psi_I\in L^2(S^m,{\bf C}).$ The BRST operator is denoted $\tilde\Omega$ and is defined by $\tilde\Omega =\widehat{\tilde \eta^a}\hat P_a.$ In local coordinates $\hat P_a$ may be identified with $\frac{1}{i}\frac{\partial}{\partial Q^a}$ as in Dirac quantization. See Woodhouse \cite{W} page 187 for details.
There exists an analogy with the DeRham complex whereby we think of the $\tilde\eta^a$'s as the coordinate one-forms $dQ^a$ and the states as complex-valued differential forms (replace all $\tilde\eta^a$'s by $dQ^a$'s). One then notes that $\tilde\Omega$ is $-idQ^a \frac {\partial}{\partial Q^a}$ which is essentially the exterior derivative. Thus the BRST cohomology will be the same as the DeRham cohomology  of $S^m$ and so is simply ${\bf C}$ at ghost numbers $-\frac{m}{2}$ and $\frac{m}{2}$ and otherwise is zero.

\bigskip
\bigskip
Finally, we consider the space of states for the combined system $M \times T^*S^m.$ We identify the extended state space with the tensor product ${\cal S}\otimes \tilde{\cal S}$ of ${\cal S}$ and $ \tilde{\cal S}$ over the space $\Lambda$ of supernumbers. Thus extended states take the form
$\Psi=\sum_{\alpha=1}^r(\psi_{\alpha}\otimes \tilde \psi_{\alpha})$ for $\psi_{\alpha}\in {\cal S}, \tilde \psi_{\alpha}\in \tilde{\cal S}.$ The ghost number operator on the extended space of states is defined by $${\cal G}^{ext}=({\cal G}\otimes 1) +(1\otimes \tilde{\cal G}).$$ 
Consequently, if $\psi \in {\cal S}^p, \tilde \psi \in \tilde{\cal S}^q,$ then 
$${\cal G}^{ext}(\psi \otimes \tilde \psi)=({\cal G}(\psi)\otimes \tilde\psi) +(\psi\otimes \tilde{\cal G}(\tilde\psi)$$$$=(-\frac{m}{2}+p)(\psi\otimes \tilde\psi)+(-\frac{m}{2}+q)(\psi\otimes\tilde\psi)=
(-m+p+q)(\psi\otimes\tilde\psi).$$ 
In particular $\psi\otimes \tilde\psi$ has ghost number zero if and only if either $p=0$ and $q=m$ or $p=m$ and $q=0.$ So states of ghost number zero lie in the space 
$$({\cal S}^0\otimes \tilde{\cal S}^m)\oplus ({\cal S}^m\otimes \tilde{\cal S}^0).$$
\bigskip
Similarly we define $\Omega^{ext}$ on the extended state space by
$$\Omega^{ext}(\Psi)=\sum_{\alpha=1}^r (\Omega(\psi_{\alpha})\otimes \tilde\psi_{\alpha})+(-1)^{\varepsilon(\psi_{\alpha})}(\psi \otimes \tilde\Omega(\tilde\psi_{\alpha}))$$
for $\Psi=\sum_{\alpha=1}^r(\psi_{\alpha}\otimes \tilde \psi_{\alpha})$ (where $\psi_{\alpha}\in {\cal S}$ and $\tilde\psi_{\alpha}\in \tilde {\cal S}$ are homogeneous). Observe that \newline
$\Omega^{ext}(\Omega^{ext}(\psi\otimes\tilde\psi))=(\Omega^2(\psi)\otimes \tilde\psi)+(-1)^{\varepsilon(\psi)}(\Omega(\psi)\otimes \tilde\Omega(\tilde\psi))$
$$+(-1)^{\varepsilon(\Omega(\psi))}(\Omega(\psi)\otimes\tilde\Omega(\tilde \psi))+
(-1)^{\varepsilon(\psi)}(-1)^{\varepsilon(\psi)}(\psi\otimes \tilde\Omega^2(\tilde\psi))=0.$$
Note that both our complexes $\{{\cal S}^p\}$ and $\{\tilde{\cal S}^q\}$ are assumed to be augmented by ${\bf C}$ where $\Omega:{\bf C}\rightarrow {\cal S}^0$ and $\tilde\Omega:{\bf C}\rightarrow \tilde{\cal S}^0$ are inclusions. Consequently, 
at ghost number zero, we see that the image of $\Omega^{ext}$ is 
$$Im(\Omega^{ext})=[{\bf C}\otimes \tilde \Omega(\tilde {\cal S}^{m-1})]
\oplus [\Omega({\cal S}^{m-1})\otimes {\bf C}]$$
and its kernel is
$$Ker(\Omega^{ext})=[ker(\Omega|_{{\cal S}^0})\otimes \tilde{\cal S}^m]\oplus
[{\cal S}^m\otimes ker(\tilde \Omega |_{\tilde{\cal S}^0})].$$

Consequently, 
$$H^0(\Omega^{ext})=\{ker(\Omega|_{{\cal S}^0})\otimes (\frac{\tilde{\cal S}^m}{ \tilde\Omega(\tilde{\cal S}^{m-1})})\}\oplus \{(\frac{{\cal S}^m}{ \Omega({\cal S}^{m-1})})\otimes ker(\tilde\Omega|_{\tilde{\cal S}^0})\}$$
$$=[H^{-\frac{m}{2}}(\Omega)\otimes H^{\frac{m}{2}}(\tilde\Omega)]\oplus 
[H^{\frac{m}{2}}(\Omega)\otimes H^{-\frac{m}{2}}(\tilde\Omega)] \quad \quad \quad \quad$$
$$=[H^{-\frac{m}{2}}(\Omega)\otimes {\bf C}]\oplus [H^{\frac{m}{2}}(\Omega)\otimes {\bf C}]=
H^{-\frac{m}{2}}(\Omega)\oplus H^{\frac{m}{2}}(\Omega).$$
\bigskip
We know from Theorem 5.4 that 
$$H^{-\frac{m}{2}}(\Omega)\cong H^{\frac{m}{2}}(\Omega)\cong {\cal PW}$$
the Hilbert space of ${\sigma}$-wave functions. Moreover, it is shown on page 320 of Henneaux and Teitelboim \cite{HT} that $\psi_0\in {\cal S}^0$ is in $ker(\Omega|_{{\cal S}^0})\cong H^{-\frac{m}{2}}(\Omega)$ if and only if $\hat G_a(\psi_0)=0.$ Thus we have the following theorem.
\begin{theorem} The BRST cohomology of $\Omega^{ext}$ at ghost number zero is two copies of the space of constrained quantum states. Moreover, this space is isomorphic to the Hilbert space of ${\sigma}$-wave functions 
${\cal PW}.$
\end{theorem}

\bigskip

To see that the constraints are correctly implemented we focus on one copy of the space of ghost zero states namely, ${\cal S}^m\otimes \tilde{\cal S}^0.$ Notice that when $\Omega^{ext} $ is restricted to this space, $\Omega^{ext}(\psi_m\otimes \tilde\psi_0)=\pm(\psi_m\otimes \tilde\Omega(\psi_0))$ on generators and $\Omega^{ext}(\psi_m\otimes \tilde\psi_0)=0$ if and only if $\tilde\Omega(\tilde \psi_0)=0$
which is the case if and only if $\tilde{\cal P}_a(\tilde\psi_0)=0$ thus implementing the choice of constraints $P_a=0$ on $T^*S^m.$ Moreover since the $\eta^b$ are the ghosts which correspond to 
the constraints $G_b=0,$ the fact that $(\hat\eta^b\otimes 1)(\psi_m\otimes \tilde\psi_0)=\hat\eta^b(\psi_m)\otimes \tilde\psi_0=0$ enforces these constraints at the quantum level on this restricted space of states. Note also that because we consider only those states in ${\cal S}^m\otimes \tilde{\cal S}^0$ for which $\Omega^{ext}=0$ we have $(1\otimes \tilde{\cal P}_a)(\psi_m\otimes \tilde\psi_0)=0.$ Finally, if we consider the special case that $M=T^*Q$ for some configuration space $Q,$ note that if $(q^i)$ are coordinates on $Q,$ then $(\hat q^i\otimes 1)(\psi_m\otimes \tilde\psi_0)=(q^i\psi_m)\otimes \tilde\psi_0)=((q^i\otimes 1)(\psi_m\otimes \tilde\psi_0)$ as one would expect in the position representation. Finally, it is clear from the theorem that the BRST cohomology on this restricted space is one of the two copies of ${\cal PW}$ described in the theorem and thus is the Hilbert space of constrained quantum states.

\bigskip
\bigskip

\noindent{\bf Remark.} Recall that the desired cohomology of $\Omega$ resides at ghost number $-\frac{m}{2}$ and that we were forced to extend $\Omega$ to $\Omega^{ext}$ in order to bring the cohomology to ghost number zero. We chose to do this by {\it essentially} tensoring $H^{-\frac{m}{2}}$
with the DeRham cohomology of a sphere. We chose a sphere because we required a space with zero DeRham cohomology except at form degree zeros and $m$ (otherwise there would be other complicating factors in the ghost zero cohomology of $\Omega^{ext}).$ The usual technique for accomplishing this is to introduce Lagrange multipliers in the Lagrangian of the theory in order to directly implement the constraints $G_a=0.$ When one has a Lagrangian this is certainly a more transparent way to implement the constraints. The Lagrange multipliers, denoted $\lambda^a$ by \cite{HT} (see page 242), become new variables which we have interpreted as coordinates $Q^a$ on the $m$-sphere. One could possibly choose some other manifold of dimension $m$ other than a sphere but the price would be paid at the cohomology level. It is not clear how one would then get rid of the unwanted states at ghost number zero.

\providecommand{\bysame}{\leavevmode\hbox to3em{\hrulefill}\thinspace}


\begin{thebibliography}{9}

 \bibitem{HT}
 M.~Henneaux and C.~Teitelboim, {\em Quantization of gauge systems},
Princeton
 Univ.~Press, 1992.

 \bibitem{W}
 N.M.J.~Woodhouse, {\em Geometric Quantization},
Oxford Univ.~Press, 1992.
 
 \bibitem{R}
A.~Rogers, { \em A global theory of supermanifolds }, J.~Math.~Phys.~
 \textbf{21} (1980), 1352--1365.



\end{thebibliography}
\end{document}